\documentclass{amsart}
\usepackage{amsmath,amssymb,amsfonts}
\usepackage[mathscr]{eucal}
\usepackage{verbatim}

\usepackage{enumerate}

\numberwithin{equation}{section}

\theoremstyle{plain}
\newtheorem{theorem}[equation]{Theorem}

\newtheorem{prop}[equation]{Proposition}
\newtheorem{lemma}[equation]{Lemma}

\theoremstyle{definition}

\newtheorem*{thank}{Acknowledgments}

\input xy
\xyoption{all}

\newcommand{\Deltaop}{{\bf \Delta}^{op}}
\newcommand{\IDeltaop}{{\bf I \Delta}^{op}}

\newcommand{\nerve}{\text{nerve}}
\newcommand{\Hom}{\text{Hom}}
\newcommand{\SSets}{\mathcal{SS}ets}

\newcommand{\LSSets}{\mathcal{LSS}ets}
\newcommand{\Tm}{\mathcal T_M}
\newcommand{\Sets}{\mathcal Sets}

\newcommand{\xu}{{\underline x}}

\newcommand{\LSSetstmi}{\mathcal {LSS}ets^{\mathcal T_{MI}}}
\newcommand{\Tmi}{\mathcal T_{MI}}

\newcommand{\Algtmi}{\mathcal Alg^{\mathcal T_{MI}}}

\newcommand{\LSSetstg}{\mathcal {LSS}ets^{\mathcal T_G}}
\newcommand{\Tg}{\mathcal T_G}

\begin{document}

\title[Erratum to ``Adding inverses"]{Erratum to ``Adding Inverses to Diagrams Encoding Algebraic Structures" and ``Adding inverses to diagrams II: Invertible homotopy theories are spaces"}

\author{Julia E. Bergner}
\email{bergnerj@member.ams.org}

\address{University of California, Riverside, CA 92521}


\thanks{The author was partially supported by NSF grant DMS-0805951.}

\maketitle

\section{Statement of previous error}

In previous work, we studied various kinds of functors $X \colon \Deltaop \rightarrow \SSets$ satisfying a Segal condition, so that the maps
\[ X_n \rightarrow \underbrace{X_1 \times_{X_0} \cdots \times_{X_0} X_1}_n \] are weak equivalences of simplicial sets for $n \geq 2$.  When we imposed the additional condition that $X_0=\Delta[0]$, such objects, called \emph{reduced Segal categories} or \emph{Segal monoids} were shown to be equivalent to simplicial monoids \cite{simpmon}.  When instead $X_0 =\amalg \Delta[0]$, some discrete simplicial set, then $X$ is a \emph{Segal category} can more generally be regarded as an up-to-homotopy model for a simplicial category with this same object set \cite{simpmon}.

These results were used in the comparison between the model structure for all Segal categories (not just with a fixed set in degree zero) and the model structure for simplicial categories \cite{thesis}.  Furthermore, the former model structure was shown to be Quillen equivalent to the model structure for complete Segal spaces, or functors $\Deltaop \rightarrow \SSets$ where the discrete level zero condition is replaced with a ``completeness" condition \cite{thesis}.

A natural question was then whether these results could be generalized to an ``invertible" and two methods were proposed in both of the papers \cite{inverses}, \cite{inverses2}.  The first was to replace $\Deltaop$ in the above definitions with a category $\IDeltaop$ in which the objects had an involution map.  However, these results were in fact incorrect, in that this involution does not adequately encode an inverse map. In this note, we clarify that this diagram should encode the structure of a monoid with involution rather  than a group, or category with involution rather than a groupoid, in the case of multiple objects.

The second approach given in these papers, given by using different projection maps as first used by Bousfield \cite{bous}, is still correct.

\begin{thank}
The author would like to thank Philip Hackney for conversations about this work.
\end{thank}

\section{Replacing groups with monoids with involution}

In the case of monoids, we consider functors $\Deltaop \rightarrow
\SSets$, where the category $\Deltaop$ has as objects finite
ordered sets $[n]=(0 \rightarrow 1 \rightarrow \cdots \rightarrow
n)$ for each $n \geq 0$ and as morphisms the opposites of the
order-preserving maps between them.  Notice that each $[n]$ can be
regarded as a category with $n+1$ objects and a single morphism $i
\rightarrow j$ whenever $i \leq j$.

In \cite{inverses} we defined a category $\IDeltaop$,
whose objects are given by small groupoids $I[n]=(0
\rightleftarrows 1 \leftrightarrows \cdots \rightleftarrows n)$
for $n \geq 0$. In other words, each $I[n]$ is a category with
$n+1$ objects and a single isomorphism between any two objects.
The morphisms of $\IDeltaop$ are generated by two sets of maps:
the opposite of the order-preserving maps from
$\Deltaop$, and an involution morphism on each $I[n]$ which
sends each $i$ to $n-i$.

The hope was that functors $\IDeltaop \rightarrow \SSets$ satisfying a Segal condition encoded a group structure.  Unfortunately, inverses are not adequately given, so such functors actually give the structure of a monoid with involution.  Thus, we change the terminology given in \cite{inverses} as follows.

In the case of ${\bf \Delta}$, the simplicial set $\Delta[n]$ is
given by the representable functor $\Hom_{\bf \Delta}(-, [n])$.
Similarly, we can define an object $I\Delta[n]$ which is
given by the representable functor $\Hom_{\bf I\Delta}(-, I[n])$.
These $n$-\emph{simplices with involution} are the standard building blocks
of the spaces we consider here.  In particular, every simplex
should be regarded as having a corresponding ``reverse" simplex.  As
with simplicial sets, we can consider the boundary of $I\Delta
[n]$, denoted $\partial I \Delta [n]$, which consists of the nondegenerate simplices of
$I \Delta [n]$ of degree less than $n$.

Thus, we can define an \emph{simplicial set with involution} to be a
functor $\IDeltaop \rightarrow \Sets$ and, more generally, an
\emph{simplicial object with involution} in a category $\mathcal C$ to
be a functor $\IDeltaop \rightarrow \mathcal C$.  We denote the
category of simplicial sets with involution by $I\SSets$.  We further
consider the case of simplicial spaces with involution, or functors
$\IDeltaop \rightarrow \SSets$. Since there is a forgetful functor
$U:I\SSets \rightarrow \SSets$ (respectively,
$U:\SSets^{I\Deltaop} \rightarrow \SSets^{\Deltaop}$), we define a
map $f$ of simplicial sets (respectively, spaces) with involution to be
a weak equivalence if $U(f)$ is a weak equivalence of simplicial
sets (respectively, spaces).

In particular, we define a \emph{Segal precategory with involution} to be an simplicial space with involution $X$ such that the simplicial set $X_0$
is discrete.  If $X_0 =\Delta[0]$, then we call it a \emph{Segal premonoid with involution}. To define a Segal category with involution, use the maps
\[ \xi_n: X_n \rightarrow \underbrace{X_1 \times_{X_0} \cdots \times_{X_0} X_1}_n \] defined in \cite[\S 4]{inverses}. Thus, a \emph{Segal category with involution} is a Segal precategory with involution $X$ such that for each $n \geq 2$
the map $\xi_n$ is a weak equivalence of simplicial sets.

Obtaining an appropriate model structure requires localization with respect to the following map:
\[ \xi_\mathcal O = \coprod_{n \geq 1}(\xi^n: \coprod_{\xu \in \mathcal O^{n+1}} (IG(n)^t_\xu \rightarrow I\Delta [n]^t_\xu)). \]

The proofs of the following two propositions continue to hold, with the necessary changes in terminology.

\begin{prop} \cite[4.1]{inverses}
There is a model category structure $\LSSets^{\IDeltaop}_{\mathcal O,f}$ on the category of
Segal precategories with involution with a fixed set $\mathcal O$ in degree zero in
which the weak equivalences and fibrations are given levelwise.
Similarly, there is a model category structure $\LSSets^{\IDeltaop}_{\mathcal O,c}$ on the
same underlying category in which the weak equivalences and
cofibrations are given levelwise.  Furthermore, we can localize
each of these model category structures with the map $\xi_\mathcal
O$ to obtain model structures $\LSSets^{\IDeltaop}_{\mathcal O,f}$ and $\LSSets^{\IDeltaop}_{\mathcal O,c}$ whose
fibrant objects are Segal categories with involution.
\end{prop}

\begin{prop} \cite[4.2]{inverses}
The adjoint pair given by the identity functor induces a Quillen
equivalence of model categories
\[ \xymatrix@1{\LSSets^{\IDeltaop}_{\mathcal O,f} \ar@<.5ex>[r] & \LSSets^{\IDeltaop}_{\mathcal O,c}. \ar@<.5ex>[l]} \]
\end{prop}

In \cite{inverses}, we claimed that there was a Quillen equivalence
\[ \LSSetstg_* \leftrightarrows \LSSets^{\IDeltaop}_{*,f}. \]  Unfortunately, we did not adequately establish that we obtained group structures using this category $\IDeltaop$.  We seek to establish that our previous proof instead gave a Quillen equivalence
\[ \LSSetstmi_* \leftrightarrows \LSSets^{\IDeltaop}_{*,f} \] where $\Tmi$ is the \emph{theory of monoids with involution}.

This theory $\Tmi$ has as objects $T_k$ which are given by the free monoid with involution on $k$ generators.  In other words, $T_k$ is free on generators $x_1, \ldots, x_k$ and $\overline{x}_1, \ldots, \overline{x}_k$ with involution $I(x_k)=\overline{x}_k$.  There is a monoid map $\tau \colon T_k \rightarrow T_k$ given by $x_j \mapsto \overline{x}_{k-j+1}$.  This map will correspond to the flip map of each object $I[k]$ in $\IDeltaop$.

Using Badzioch's theorem from \cite{bad}, such a Quillen equivalence will complete the proof of the following theorem.

\begin{theorem} \label{main}
The model category structure $\Algtmi$ is Quillen equivalent to the
model category structure $\LSSets^{\IDeltaop}_{*,f}$.
\end{theorem}

As in \cite{simpmon}, we prove this theorem using several lemmas.
Note that in the model structure $\LSSets^{\IDeltaop}_{*,c}$, we denote by
$L_1$ the localization, or fibrant replacement functor.
Analogously, we denote by $L_2$ the localization functor in
$\LSSetstmi_*$.

The first step in the proof of the theorem is to show what the
localization functor $L_1$ does to the $n$-simplex with involution
$I\Delta[n]^t$.  By $I\nerve(-)^t$, we denote the representable
functor $\Hom(I[n],-)$, viewed as a transposed constant simplicial space.  It is here, Lemma \cite[4.5]{inverses}, that the major error occurred.  The correct statement is as follows.

\begin{prop}  \label{nerve}
Let $F_n$ denote the free monoid with involution on $n$ generators.  Then in
$\LSSets^{\IDeltaop}_{*,c}$, $L_1 I\Delta [n]^t_*$ is weakly equivalent to $I\nerve (F_n)^t$ for each $n \geq 0$.
\end{prop}

In the proof, we defined a filtration of $I\nerve(F_1)^t$ as follows:
\[ \Psi_k (I\nerve(F_1)^t)_j=\left\{(x^{n_1}|\cdots |x^{n_j}) \mid \sum_{\ell=1}^j |n_\ell| \leq k \right\} \]
where $x$ and its ``inverse" $x^{-1}$ denote the two nondegenerate 1-simplices of $I \Delta[1]^t_*= \Psi_1$.  The problem is that we assumed here that $x$ could be canceled with $x^{-1}$, when there is no structure built in to $\IDeltaop$ to make this cancelation possible.  In short, we do not really have the structure of a group, but only of a monoid with involution.  Therefore, we instead define the filtration of $I\nerve(F_1)^t$ by
\[ \Psi_k (I\nerve(F_1)^t)_j=\left\{(w_1|\cdots |w_j) \mid \ell(w_1 \cdots w_j) \leq k \right\} \] where the $w_i$ are words in $x$ and $\overline{x}$, and $\ell$ denotes word length.

With this modification, the previous proof goes through as before.

To obtain the appropriate adjoint pair for the Quillen equivalence, we first define a functor $J \colon \IDeltaop \rightarrow \Tmi$.  On objects, this functor is defined by $I{n} \mapsto T_n$.  On face and degeneracy maps coming from $\Deltaop$, this functor behaves the same as the functor $\Deltaop \rightarrow \Tm$ as defined in \cite{simpmon}.  Specifically, the coface and codegeneracy maps maps in ${\bf I\Delta}$ are given by
\[ d^i(x_k)=
\begin{cases}
x_k & i>k \\
x_kx_{k+1} & i=k \\
x_{k+1} & i<k
\end{cases} \text{ and }
s^i(x_k)=
\begin{cases}
x_k & i \leq k \\
e & i=k-1 \\
x_{k-1} & i<k-1.
\end{cases} \]
The involution map is sent to the map $\tau \colon T_\ell \rightarrow T_\ell$ defined by $x_k \mapsto y_{\ell-k+1}$.

This functor $J$ induces a map
\[ J^\ast \colon \LSSetstmi_{\ast,f} \rightarrow \LSSets^{\IDeltaop}_{\ast,f} \] for which we have a left adjoint $J_\ast$ via left Kan extension.

As in \cite{inverses}, define $IM[k]$ to be the functor $\Tmi \rightarrow \SSets$ given by
\[ F_n \mapsto \Hom_{\Tmi}(F_k, F_n)=(F_k)^n, \] and let $H = I\nerve(F_k)^t$.  The following results continue to hold, replacing $\Tg$ with $\Tmi$.

\begin{lemma} \cite[4.9]{inverses}
In $\LSSets^{\Tmi}_\ast$, $L_2J_*(H)$ is weakly equivalent to $IM[k]$.
\end{lemma}

\begin{prop} \cite[4.10]{inverses}
For any object $X$ in $\LSSets^{\IDeltaop}_{\ast, c}$, $L_1X$ is weakly equivalent to $J^*L_2J_*X$.
\end{prop}

\begin{prop} \cite[4.11, 4.3]{inverses}
The adjoint pair
\[ \xymatrix@1{J_*: \LSSets^{\IDeltaop}_{*,f} \ar@<.5ex>[r] & \SSets^{\Tmi}_* :J^* \ar@<.5ex>[l]} \]
is a Quillen equivalence.
\end{prop}

The more general fixed-object case, as well as the more general case as developed in \cite{inverses2}, can be corrected by considering categories with involution rather than groupoids whenever $\IDeltaop$ is the indexing category.

\end{document}